\documentclass [12pt]{article}
\usepackage{graphicx}
\begin{document}

\title{Frequency Response of Uncertain Systems: Strong Kharitonov-Like Results\thanks{Supported by National Natural Science Foundation (69925307) and National Key Project of China.} }
\author{Long Wang\\{\small Center for Systems and Control, Department of Mechanics and Engineering Science}\\{\small Peking University, Beijing 100871, China}}
\date{}
\maketitle
%\baselineskip 24pt
%\large

\small

\begin{abstract}

In this paper, we study the frequency response of uncertain systems using Kharitonov stability theory on {\it first order} complex polynomial set. For an interval transfer function, we show that the minimal real part of the frequency response at any fixed frequency is attained at some prescribed vertex transfer functions. By further geometric and algebraic analysis, we identify an index for strict positive realness of interval transfer functions. Some extensions and applications in positivity verification and robust absolute stability of feedback control systems are also presented.


\end{abstract}

\vspace{0.6cm}

\section{Introduction}

Motivated by the seminal theorem of Kharitonov on robust stability of interval
polynomials\cite{ack, khar}, a number of papers on robustness analysis of uncertain
systems have been published in the past few years\cite{bart, tsyp, Barmish, wang3,
Ack1, Ack2, chap, das, Ran, Holl, wang, Wang1}. Kharitonov's theorem states that the Hurwitz stability of
a real (or complex) interval polynomial family can be guaranteed by the Hurwitz
stability of four (or eight) prescribed critical vertex polynomials in this family.
This result is significant since it reduces checking stability of infinitely many
polynomials to checking stability of finitely many polynomials, and the number of
critical vertex polynomials need to be checked is independent of the order of the
polynomial family. An important extension of Kharitonov's theorem is the edge theorem
discovered by Bartlett, Hollot and Huang\cite{bart}. The edge theorem states that
the stability of a polytope of polynomials can be guaranteed by the stability of its
one-dimensional exposed edge polynomials. The significance of the edge theorem is
that it allows some (affine) dependency among polynomial coefficients, and applies to
more general stability regions, e.g., unit circle, left sector, shifted half plane,
hyperbola region, etc. When the dependency among polynomial coefficients is nonlinear,
however, Ackermann shows that checking a subset of a polynomial family generally can
not guarantee the stability of the entire family\cite{Ack1, Ack2}.

In this paper, we study the frequency response of uncertain systems using Kharitonov stability theory on {\it first order} complex polynomial set. For an interval transfer function, we show that the minimal real part of the frequency response at any fixed frequency is attained at some prescribed vertex transfer functions. By further geometric and algebraic analysis, we identify an index for strict positive realness of interval transfer functions. Some extensions and applications in positivity verification and robust absolute stability of feedback control systems are also presented.

\section{Preliminaries}

Denote the $m$-th, $n$-th $(m \leq n)$ order real interval polynomials $K^m (s), K^n (s)$ as

\begin{equation}
K^m (s) = \{ g(s) | g(s) = \Sigma_{i=0}^{m} a_i s^i, a_i \in [ a_i^- , a_i^+ ], i=0,1,......,m \}
\end{equation}

\begin{equation}
K^n (s) = \{ f(s) | f(s) = \Sigma_{j=0}^{n} b_j s^j, b_j \in [ b_j^- , b_j^+ ], j=0,1,......,n \}
\end{equation}
For any $g(s) \in K^m (s)$, we have

\begin{equation}
g(s) = \alpha_g (s^2) + s \beta_g (s^2)
\end{equation}
where

\begin{equation}
\alpha_g (s^2) = a_0 + a_2 s^2 + a_4 s^4 + a_6 s^6 + ......
\end{equation}

\begin{equation}
\beta_g (s^2) = a_1 + a_3 s^2 + a_5 s^4 + a_7 s^6 + ......
\end{equation}

For the interval polynomial $K^m (s)$, denote

\begin{equation}
\alpha_g^{(1)} (s^2) = a_0^- + a_2^+ s^2 + a_4^- s^4 + a_6^+ s^6 + ......
\end{equation}

\begin{equation}
\alpha_g^{(2)} (s^2) = a_0^+ + a_2^- s^2 + a_4^+ s^4 + a_6^- s^6 + ......
\end{equation}

\begin{equation}
\beta_g^{(1)} (s^2) = a_1^- + a_3^+ s^2 + a_5^- s^4 + a_7^+ s^6 + ......
\end{equation}

\begin{equation}
\beta_g^{(2)} (s^2) = a_1^+ + a_3^- s^2 + a_5^+ s^4 + a_7^- s^6 + ......
\end{equation}
and denote $g_{ij} (s) \in K^m (s)$ as

\begin{equation}
g_{ij}(s) = \alpha_g^{(i)} (s^2) + s \beta_g^{(j)} (s^2), \quad i,j = 1,2
\end{equation}

For the interval polynomial $K^n (s)$, the corresponding $\alpha_f^{(i)} (s^2)$, $\beta_f^{(j)} (s^2)$ and $f_{ij} (s) \in K^n (s)$ can be defined analogously.

Denote the set of all Hurwitz stable polynomials as $H$.

A transfer function $\frac{g(s)}{f(s)}$ is said to be strictly positive real (SPR), denoted as $\frac{g(s)}{f(s)} \in SPR$, if

$1) \quad f(s) \in H$

$2) \quad \Re \frac{g(j \omega)}{f(j \omega)} > 0, \forall \omega \in R$

A transfer function set $J$ is said to be SPR, denoted as $J \in SPR$, if every member in $J$ is SPR.

\vspace{0.6cm}

\noindent {\bf Lemma 1}\cite{ack}

For any fixed $\omega \in R$ and any $g(s) \in K^m (s)$, $f(s) \in K^n (s)$

\begin{equation}
\alpha_g^{(1)} (- \omega^2) \leq \alpha_g (- \omega^2) \leq \alpha_g^{(2)} (- \omega^2)
\end{equation}

\begin{equation}
\alpha_f^{(1)} (- \omega^2) \leq \alpha_f (- \omega^2) \leq \alpha_f^{(2)} (- \omega^2)
\end{equation}

\begin{equation}
\beta_g^{(1)} (- \omega^2) \leq \beta_g (- \omega^2) \leq \beta_g^{(2)} (- \omega^2)
\end{equation}

\begin{equation}
\beta_f^{(1)} (- \omega^2) \leq \beta_f (- \omega^2) \leq \beta_f^{(2)} (- \omega^2)
\end{equation}

\vspace{0.6cm}

\noindent {\bf Lemma 2}

For any fixed $\omega, \beta \in R$, if $f(j \omega) \neq 0$, then

\begin{equation}
\Re \frac{g(j \omega) - \beta f(j \omega)}{f(j \omega)} > 0
\end{equation}
if and only if

\begin{equation}
f(j \omega) s + g(j \omega) - \beta f(j \omega) \in H
\end{equation}

Proof: For any fixed $\omega, \beta \in R$, $g(j \omega) - \beta f(j \omega)$ and $f(j \omega)$ are fixed complex numbers. Thus $f(j \omega) s + [g(j \omega) - \beta f(j \omega)]$ is a first order complex polynomial, and

\begin{equation}
\begin{array}{ll}
& f(j \omega) s + g(j \omega) - \beta f(j \omega) \in H\\[6mm]
\Longleftrightarrow & \Re \frac{-[g(j \omega) - \beta f(j \omega)]}{f(j \omega)} < 0\\[6mm]
\Longleftrightarrow & \Re \frac{g(j \omega) - \beta f(j \omega)}{f(j \omega)} > 0
\end{array}
\end{equation}
This completes the proof.

Consider the interval complex numbers $c_0 + jd_0, c_1 + jd_1$, where $c_i \in [c_i^-, c_i^+], d_i \in [d_i^-, d_i^+], i=1,2$. Define the sign functions

\begin{equation}
sgn[c_i] = \left\{ \begin{array}{ll}
                   1 & \quad \quad c_i=c_i^-\\[6mm]
                   2 & \quad \quad c_i=c_i^+
                   \end{array}
           \right.
\end{equation}

\begin{equation}
sgn[d_i] = \left\{ \begin{array}{ll}
                   1 & \quad \quad d_i=d_i^-\\[6mm]
                   2 & \quad \quad d_i=d_i^+
                   \end{array}
           \right.
\end{equation}
and define the index sets $I_1, I_2, I_3$ as

\begin{equation}
\begin{array}{l}
I_1 = \{ (1222), (1221), (2221), (2211), (2111), (2112), (1112), (1122), (1211), (2212), (2122), (1121) \}\\[6mm]
I_2 = \{ (1112), (1222), (2111), (2221), (1121), (1211), (2122), (2212) \}\\[6mm]
I_3 = \{ (1111), (1212), (2222), (2121), (1112), (1222), (2221), (2111), (1211), (2212), (2122), (1121) \}
\end{array}
\end{equation}

For any two polynomials $h^{(1)} (s), h^{(2)} (s)$, denote their convex combination as

\begin{equation}
L[h^{(1)}, h^{(2)}] = \{ \lambda h^{(1)} (s) + (1- \lambda)  h^{(2)} (s) | \lambda \in [0, 1] \}
\end{equation}

\noindent {\bf Lemma 3}

For any fixed $\beta > 0$, the first order complex polynomial set

\begin{equation}
W_1 (s) := \{ (c_1 + jd_1) (s- \beta) + (c_0 + jd_0) | c_i \in [c_i^-, c_i^+], d_i \in [d_i^-, d_i^+], i=1,2 \} \subset H
\end{equation}
if and only if

\begin{equation}
\{ (c_1 + jd_1) (s- \beta) + (c_0 + jd_0) | (sgn[c_0] \quad sgn[d_0] \quad sgn[c_1] \quad sgn[d_1]) \in I_1 \} \subset H
\end{equation}

Proof: Necessity is obvious, to prove sufficiency, by Zero Exclusion Principle\cite{ack}, we only need to prove

\begin{equation}
0 \not\in W_1 (j \omega), \quad \forall \omega \in R
\end{equation}
Since $W_1 (s)$ has a fixed order, for sufficiently large $\omega_{\infty}$

\begin{equation}
0\not\in W_1 (j \omega_{\infty})
\end{equation}
Hence, by continuity, we only need to prove

\begin{equation}
0 \not\in \partial W_1 (j \omega), \quad \forall \omega \in R
\end{equation}
where $\partial W_1 (j \omega)$ stands for the boundary set of $W_1 (j \omega)$ in the complex plane.

Apparently, the interval complex numbers $c_0 + jd_0, c_1 + jd_1$ are two rectangles in the complex plane, with edges parallel to the coordinate axes (Fig. 1).
When $\omega \geq 0$, $(j \omega - \beta)(c_1 + jd_1)$ is a rotated rectangle (Fig. 2). Hence, $W_1 (j \omega)$ is a octagon, produced by addition of two rectangles in the complex plane (Fig. 3).

To prove $0 \not\in \partial W_1 (j \omega), \forall \omega \geq 0$, suppose on the contrary that there exists $\omega_0 \geq 0$ such that $0 \in \partial W_1 (j \omega_0)$. Let

\begin{equation}
h_1 (s) = (c_1^+ + jd_1^+) (s- \beta) + (c_0^- + jd_0^+)
\end{equation}

\begin{equation}
h_2 (s) = (c_1^+ + jd_1^-) (s- \beta) + (c_0^- + jd_0^+)
\end{equation}
Without loss of generality, suppose

\begin{equation}
0 \in L[h_1, h_2] (j \omega_0)
\end{equation}
Namely, there exists $\lambda_0 \in (0, 1)$ such that

\begin{equation}
\lambda_0 h_1 (j \omega_0) + (1- \lambda_0)  h_2 (j \omega_0) = 0
\end{equation}
Since $h_1 (s), h_2 (s) \in H$, we have

\begin{equation}
\frac{d}{d \omega} \arg h_i (j \omega) > 0, \quad \forall \omega \in R, \quad i = 1,2
\end{equation}
Thus, we have

\begin{equation}
\begin{array}{ll}
  & \frac{d}{d \omega} \arg [h_2 - h_1] (j \omega) |_{\omega_0}\\[6mm]
= & (1- \lambda_0) \frac{d}{d \omega} \arg h_1 (j \omega) |_{\omega_0} + \lambda_0 \frac{d}{d \omega} \arg h_2 (j \omega) |_{\omega_0} > 0
\end{array}
\end{equation}
On the other hand, since

\begin{equation}
h_2 (s) - h_1 (s) = j (d_1^- - d_1^+) (s- \beta)
\end{equation}
and $\beta > 0$, obviously we have

\begin{equation}
\frac{d}{d \omega} \arg [h_2 - h_1] (j \omega) < 0, \quad \forall \omega \in R
\end{equation}
Hence, a contradiction arises, which means that $0 \in \partial W_1 (j \omega_0)$ is impossible.

\noindent When $\omega < 0$, $W_1 (j \omega)$ is still a octagon in the complex plane (Fig. 4). However, four of the eight vertices are different from the previous ones. By similar analysis, we have

\begin{equation}
0 \not\in \partial W_1 (j \omega), \quad \forall \omega < 0
\end{equation}
This completes the proof.

The following result is a direct consequence of Kharitonov Theorem for interval complex polynomials.

\vspace{0.6cm}

\noindent {\bf Lemma 4}\cite{ack, khar}

The first order interval complex polynomial set

\begin{equation}
W_2 (s) := \{ (c_1 + jd_1) s + (c_0 + jd_0) | c_i \in [c_i^-, c_i^+], d_i \in [d_i^-, d_i^+], i=1,2 \} \subset H
\end{equation}
if and only if

\begin{equation}
\{ (c_1 + jd_1) s + (c_0 + jd_0) | (sgn[c_0] \quad sgn[d_0] \quad sgn[c_1] \quad sgn[d_1]) \in I_2 \} \subset H
\end{equation}

\section{Pointwise Strict Positive Realness}

\vspace{0.6cm}

\noindent {\bf Theorem 1}

For any fixed $\omega \in R$, if $0 \not\in K^n (j \omega)$ and

\begin{equation}
\min \{ \Re \frac{g_{i_1 j_1} (j \omega)}{f_{i_2 j_2} (j \omega)} | (i_1 \quad j_1 \quad i_2 \quad j_2) \in I_1 \} := \beta_0 > 0
\end{equation}
Then

\begin{equation}
\min \{ \Re \frac{g (j \omega)}{f (j \omega)} | g(s) \in K^m (s), f(s) \in K^n (s) \} = \beta_0
\end{equation}

Proof: Since $g_{ij} (s) \in K^m (s), f_{ij} (s) \in K^n (s), i,j=1,2$, we have

\begin{equation}
\min \{ \Re \frac{g (j \omega)}{f (j \omega)} | g(s) \in K^m (s), f(s) \in K^n (s) \} \leq \beta_0
\end{equation}
Suppose

\begin{equation}
\min \{ \Re \frac{g (j \omega)}{f (j \omega)} | g(s) \in K^m (s), f(s) \in K^n (s) \} := \beta_1 < \beta_0
\end{equation}
Since $\beta_0 > 0$, there exists $\beta_2 > 0$ such that $\beta_1 < \beta_2 < \beta_0$. Hence, for any $(i_1 \quad j_1 \quad i_2 \quad j_2) \in I_1$, we have

\begin{equation}
\Re \frac{g_{i_1 j_1} (j \omega)}{f_{i_2 j_2} (j \omega)} \geq \beta_0 > \beta_2 > 0
\end{equation}
Namely

\begin{equation}
\Re \frac{g_{i_1 j_1} (j \omega) - \beta_2 f_{i_2 j_2} (j \omega)}{f_{i_2 j_2} (j \omega)} > 0
\end{equation}
By Lemma 2, we have

\begin{equation}
f_{i_2 j_2} (j \omega) s + g_{i_1 j_1} (j \omega) - \beta_2 f_{i_2 j_2} (j \omega) \in H, \quad \forall (i_1 \quad j_1 \quad i_2 \quad j_2) \in I_1
\end{equation}

Consider the first order complex polynomial set

\begin{equation}
W_3 (s) := \{ f (j \omega) s + g (j \omega) - \beta_2 f (j \omega) | g(s) \in K^m (s), f(s) \in K^n (s) \}
\end{equation}
By Lemma 1, when $\omega \geq 0$, we have

\begin{equation}
\alpha_f^{(1)} (- \omega^2) \leq \Re f (j \omega) \leq \alpha_f^{(2)} (- \omega^2)
\end{equation}

\begin{equation}
\omega \beta_f^{(1)} (- \omega^2) \leq \Im f (j \omega) \leq \omega \beta_f^{(2)} (- \omega^2)
\end{equation}

\begin{equation}
\alpha_g^{(1)} (- \omega^2) \leq \Re g (j \omega) \leq \alpha_g^{(2)} (- \omega^2)
\end{equation}

\begin{equation}
\omega \beta_g^{(1)} (- \omega^2) \leq \Im g (j \omega) \leq \omega \beta_g^{(2)} (- \omega^2)
\end{equation}
By Lemma 3, $W_3 (s) \subset H$. When $\omega < 0$, the two inequalities on the imaginary parts of $f(j \omega), g(j \omega)$ above will be reversed. By Lemma 3, $W_3 (s) \subset H$. Hence, for any fixed $\omega \in R, f(s) \in K^n (s), g(s) \in K^m (s)$, we have

\begin{equation}
f (j \omega) s + g (j \omega) - \beta_2 f (j \omega) \in H
\end{equation}
By Lemma 2, we have

\begin{equation}
\Re \frac{g (j \omega) - \beta_2 f (j \omega)}{f (j \omega)} > 0, \quad \forall f(s) \in K^n (s), g(s) \in K^m (s)
\end{equation}
Namely

\begin{equation}
\Re \frac{g (j \omega)}{f (j \omega)} > \beta_2, \quad \forall f(s) \in K^n (s), g(s) \in K^m (s)
\end{equation}
Namely

\begin{equation}
\min \{ \Re \frac{g (j \omega)}{f (j \omega)} | g(s) \in K^m (s), f(s) \in K^n (s) \} = \beta_1 > \beta_2
\end{equation}
which contradicts $\beta_1 < \beta_2 < \beta_0$. This completes the proof.

\vspace{0.6cm}

\noindent {\bf Corollary 1a}

If $f_{ij} (s) \in H, i,j=1,2$ and

\begin{equation}
\min \{ \inf_{\omega \in R} \Re \frac{g_{i_1 j_1} (j \omega)}{f_{i_2 j_2} (j \omega)} | (i_1 \quad j_1 \quad i_2 \quad j_2) \in I_1 \} := \gamma_0 > 0
\end{equation}
Then

\begin{equation}
\min \{ \inf_{\omega \in R} \Re \frac{g (j \omega)}{f (j \omega)} | g(s) \in K^m (s), f(s) \in K^n (s) \} = \gamma_0
\end{equation}

Proof: Since $f_{ij} (s) \in H, i,j=1,2$, by Kharitonov Theorem\cite{ack, khar}, $K^n (s) \subset H$. Hence

\begin{equation}
0 \not\in K^n (j \omega), \quad \forall \omega \in R
\end{equation}
Moreover, since $g_{ij} (s) \in K^m (s), f_{ij} (s) \in K^n (s), i,j=1,2$, we have

\begin{equation}
\min \{ \inf_{\omega \in R} \Re \frac{g (j \omega)}{f (j \omega)} | g(s) \in K^m (s), f(s) \in K^n (s) \} := \gamma_1 \leq \gamma_0
\end{equation}
Suppose $\gamma_1 < \gamma_0$, since $\gamma_0 > 0$, there exists $\gamma_2 > 0$ such that $\gamma_1 < \gamma_2 < \gamma_0$. Since $\gamma_0 > \gamma_2 > 0$, for any fixed $\omega \in R$, we have

\begin{equation}
\Re \frac{g_{i_1 j_1} (j \omega)}{f_{i_2 j_2} (j \omega)} > \gamma_2 > 0, \quad \forall (i_1 \quad j_1 \quad i_2 \quad j_2) \in I_1
\end{equation}
By Theorem 1, we have

\begin{equation}
\Re \frac{g (j \omega)}{f (j \omega)} > \gamma_2 > 0, \quad \forall f(s) \in K^n (s), g(s) \in K^m (s)
\end{equation}
Hence, we have

\begin{equation}
\inf_{\omega \in R} \Re \frac{g (j \omega)}{f (j \omega)} \geq \gamma_2
\end{equation}
Namely

\begin{equation}
\min \{ \inf_{\omega \in R} \Re \frac{g (j \omega)}{f (j \omega)} | g(s) \in K^m (s), f(s) \in K^n (s) \} = \gamma_1 \geq \gamma_2
\end{equation}
which contradicts $\gamma_1 < \gamma_2 < \gamma_0$. This completes the proof.

By similar analysis, we have

\vspace{0.6cm}

\noindent {\bf Corollary 1b}

If $\forall \omega \in [\omega_1, \omega_2], 0 \not\in K^n (j \omega)$ and

\begin{equation}
\min \{ \inf_{\omega \in [\omega_1, \omega_2]} \Re \frac{g_{i_1 j_1} (j \omega)}{f_{i_2 j_2} (j \omega)} | (i_1 \quad j_1 \quad i_2 \quad j_2) \in I_1 \} := \gamma_0 > 0
\end{equation}
Then

\begin{equation}
\min \{ \inf_{\omega \in [\omega_1, \omega_2]} \Re \frac{g (j \omega)}{f (j \omega)} | g(s) \in K^m (s), f(s) \in K^n (s) \} = \gamma_0
\end{equation}

\vspace{0.6cm}

\noindent {\bf Theorem 2}

For any fixed $\omega \in R$, if $0 \not\in K^n (j \omega)$, then

\begin{equation}
\min \{ \Re \frac{g (j \omega)}{f (j \omega)} | g(s) \in K^m (s), f(s) \in K^n (s) \} > 0
\end{equation}
if and only if

\begin{equation}
\min \{ \Re \frac{g_{i_1 j_1} (j \omega)}{f_{i_2 j_2} (j \omega)} | (i_1 \quad j_1 \quad i_2 \quad j_2) \in I_2 \} > 0
\end{equation}

Proof: Necessity: Obvious.

Sufficiency: Since

\begin{equation}
\min \{ \Re \frac{g_{i_1 j_1} (j \omega)}{f_{i_2 j_2} (j \omega)} | (i_1 \quad j_1 \quad i_2 \quad j_2) \in I_2 \} > 0
\end{equation}
By Lemma 2, for any fixed $\omega \in R$

\begin{equation}
f_{i_2 j_2} (j \omega) s + g_{i_1 j_1} (j \omega) \in H, \quad \forall (i_1 \quad j_1 \quad i_2 \quad j_2) \in I_2
\end{equation}

Consider the first order interval complex polynomial set

\begin{equation}
W_4 (s) := \{ f (j \omega) s + g (j \omega) | g(s) \in K^m (s), f(s) \in K^n (s) \}
\end{equation}
By the proof of Theorem 1 and by Lemma 4, $W_4 (s) \subset H$. Hence, for any fixed $\omega \in R, f(s) \in K^n (s), g(s) \in K^m (s)$, we have

\begin{equation}
f (j \omega) s + g (j \omega) \in H
\end{equation}
By Lemma 2, we have

\begin{equation}
\Re \frac{g (j \omega)}{f (j \omega)} > 0, \quad \forall f(s) \in K^n (s), g(s) \in K^m (s)
\end{equation}
This completes the proof.

\vspace{0.6cm}

\noindent {\bf Corollary 2}\cite{wang3, chap}

\begin{equation}
\{ \frac{g(s)}{f(s)} | g(s) \in K^m (s), f(s) \in K^n (s) \} \subset SPR
\end{equation}
if and only if

\begin{equation}
\{ \frac{g_{i_1 j_1} (s)}{f_{i_2 j_2} (s)} | (i_1 \quad j_1 \quad i_2 \quad j_2) \in I_2 \} \subset SPR
\end{equation}

Proof: Necessity: Obvious.

Sufficiency: Since

\begin{equation}
\{ \frac{g_{i_1 j_1} (s)}{f_{i_2 j_2} (s)} | (i_1 \quad j_1 \quad i_2 \quad j_2) \in I_2 \} \subset SPR
\end{equation}
We have

\begin{equation}
f_{i_2 j_2} (s) \in H, \quad i_2, j_2 = 1, 2
\end{equation}
and

\begin{equation}
\min \{ \Re \frac{g_{i_1 j_1} (j \omega)}{f_{i_2 j_2} (j \omega)} | (i_1 \quad j_1 \quad i_2 \quad j_2) \in I_2 \} > 0
\end{equation}
By Kharitonov Theorem\cite{ack, khar}, $K^n (s) \subset H$. Hence, $0 \not\in K^n (j \omega)$. By Theorem 2, we have

\begin{equation}
\min \{ \Re \frac{g (j \omega)}{f (j \omega)} | g(s) \in K^m (s), f(s) \in K^n (s) \} > 0
\end{equation}
This completes the proof.

Theorem 2 is stronger than Corollary 2, since Theorem 2 reveals a pointwise property of the frequency response. This can be illustrated in the following example.

\vspace{0.6cm}

\noindent {\bf Example 1}

Consider the interval transfer function

\begin{equation}
\frac{[3, 5]s + [-7, 9]}{[5, 8]s +[1, 2]}
\end{equation}
and suppose $\omega = 1$. Then it is easy to verify that all of the eight vertex transfer functions in Theorem 2 have positive real parts at this frequency. Hence by Theorem 2, all of the transfer functions in this interval family have positive real parts at this frequency. However, Corollary 2 does not apply in this case, since some transfer functions in this family, like $\frac{4s-6}{7s+2}$, are not strictly positive real. Similar results can be shown for the following interval transfer functions

\begin{equation}
\frac{[-2, 5]s^2 + [3, 5]s + [-2, 7]}{[2, 4]s^5 + [3, 4]s +[1, 2]}
\end{equation}

\begin{equation}
\frac{[2, 3]s^9 + [1, 2]s + [-7, 9]}{[-3, -2]s^3 + [3, 5]s +[1, 2]}
\end{equation}

\section{Further Extensions and Applications}

\subsection{Index of Strict Positive Realness}

\vspace{0.6cm}

\noindent {\bf Lemma 5}\cite{ack}

Consider the real polynomial

\begin{equation}
h(s, \lambda) = h^{(0)} (s) + \lambda (\alpha s + \beta), \quad \lambda \in [0, 1]
\end{equation}
with constant order. Then

\begin{equation}
h(s, \lambda) \in H, \quad \forall \lambda \in [0, 1]
\end{equation}
if and only if

\begin{equation}
h(s, 0), h(s, 1) \in H
\end{equation}

\vspace{0.6cm}

\noindent {\bf Lemma 6}

For any fixed $\beta > 0$, the first order real polynomial set

\begin{equation}
W_5 (s) := \{ c_1 (s+ \beta) + c_0 | c_i \in [c_i^-, c_i^+], i=1,2 \} \subset H
\end{equation}
if and only if

\begin{equation}
\{ c_1 (s+ \beta) + c_0 | c_i \in \{ c_i^-, c_i^+ \}, i=1,2 \} \subset H
\end{equation}

Proof: Necessity: Obvious.

Sufficiency: Suppose $c_1 \in [c_1^-, c_1^+]$ is fixed, if

\begin{equation}
c_1 (s+ \beta) + c_0^- \in H
\end{equation}

\begin{equation}
c_1 (s+ \beta) + c_0^+ \in H
\end{equation}
Then, by Lemma 5, we have

\begin{equation}
c_1 (s+ \beta) + c_0 \in H, \quad \forall c_0 \in [c_0^-, c_0^+]
\end{equation}
Suppose $c_0 \in [c_0^-, c_0^+]$ is fixed, if

\begin{equation}
c_1^- (s+ \beta) + c_0 \in H
\end{equation}

\begin{equation}
c_1^+ (s+ \beta) + c_0 \in H
\end{equation}
Then, by Lemma 5, we have

\begin{equation}
c_1 (s+ \beta) + c_0 \in H, \quad \forall c_1 \in [c_1^-, c_1^+]
\end{equation}
This completes the proof.

\vspace{0.6cm}

\noindent {\bf Lemma 7}\cite{chap}

Suppose $f(s) \in H$, then $\frac{g(s)}{f(s)} \in SPR$ if and only if

\begin{equation}
\begin{array}{ll}
1) & \Re \frac{g(0)}{f(0)} > 0\\[6mm]
2) & g(s) \in H\\[6mm]
3) & f(s) + j \alpha g(s) \in H, \quad \forall \alpha \in R
\end{array}
\end{equation}

\vspace{0.6cm}

\noindent {\bf Lemma 8}

For any fixed $\beta > 0, \gamma \in R \backslash \{ 0 \}$, the complex polynomial set

\begin{equation}
W_6 (s) := \{ g(s) + (\beta + j \gamma) f(s) | g(s) \in K^m (s), f(s) \in K^n (s) \} \subset H
\end{equation}
if and only if

\begin{equation}
\{ g_{i_1 j_1} (s) + (\beta + j \gamma) f_{i_2 j_2}(s) | (i_1 \quad j_1 \quad i_2 \quad j_2) \in I_3 \} \subset H
\end{equation}

Proof: Necessity: Obvious.

Sufficiency: When $\gamma > 0, \arg (\beta + j \gamma) \in (0,
\frac{\pi}{2})$. When $\omega \geq 0$, by Lemma 1, $K^m (j
\omega), K^n (j \omega)$ are rectangles in the complex plane,
with edges parallel to the coordinate axes (Fig. 5). The four

vertices of $K^m (j \omega)$ are $g_{11} (j \omega), g_{12} (j
\omega), g_{22} (j \omega), g_{21} (j \omega)$. $(\beta + j
\gamma) K^n (j \omega)$ is produced by rotating the rectangle
$K^n (j \omega)$ counterclockwisely by $\arg (\beta + j \gamma)$
with respect to the origin of the coordinate axes, and then
scaling by $|\beta + j \gamma|$ (Fig. 6). Hence, $W_6 (j \omega)
= K^m (j \omega) + (\beta + j \gamma) K^n (j \omega)$ is a
octagon in the complex plane (Fig. 7), with edges parallel to
either the edges of $K^m (j \omega)$ or the edges of $(\beta + j
\gamma) K^n (j \omega)$. Hence, the slopes of the edges of $W_6
(j \omega)$ are invariant with respect to $\omega$. Moreover,

since $\{ g_{i_1 j_1} (s) + (\beta + j \gamma) f_{i_2 j_2}(s) |
(i_1 \quad j_1 \quad i_2 \quad j_2) \in I_3 \} \subset H$, by
similar analysis as in the proof of Lemma 3, we have

\begin{equation}
0 \not\in \partial W_6 (j \omega), \quad \forall \omega \geq 0
\end{equation}
The other three cases ($\gamma > 0, \omega < 0; \gamma < 0, \omega \geq 0; \gamma < 0, \omega < 0$) can be analyzed analogously. Henceforth, for any fixed $\omega \in R, \beta > 0, \gamma \in R \backslash \{ 0 \}$, we have

\begin{equation}
0 \not\in \partial W_6 (j \omega)
\end{equation}
Moreover, since $\gamma \neq 0$, $W_6 (s)$ has a constant order. Hence

\begin{equation}
0 \not\in W_6 (j \omega)
\end{equation}
By Zero Exclusion Principle\cite{ack}, $W_6 (s) \subset H$. This completes the proof.

\vspace{0.6cm}

\noindent {\bf Theorem 3}

If $f_{ij} (s) \in H, i,j=1,2$ and

\begin{equation}
\min \{ \inf_{\omega \in R} \Re \frac{g_{i_1 j_1} (j \omega)}{f_{i_2 j_2} (j \omega)} | (i_1 \quad j_1 \quad i_2 \quad j_2) \in I_3 \} := \gamma_0 < 0
\end{equation}
Then

\begin{equation}
\min \{ \inf_{\omega \in R} \Re \frac{g (j \omega)}{f (j \omega)} | g(s) \in K^m (s), f(s) \in K^n (s) \} = \gamma_0
\end{equation}

Proof: Since $f_{ij} (s) \in H, i,j=1,2$, by Kharitonov Theorem\cite{ack, khar}, $K^n (s) \subset H$. Moreover, since $g_{ij} (s) \in K^m (s), f_{ij} (s) \in K^n (s), i,j=1,2$, we have

\begin{equation}
\min \{ \inf_{\omega \in R} \Re \frac{g (j \omega)}{f (j \omega)} | g(s) \in K^m (s), f(s) \in K^n (s) \} := \gamma_1 \leq \gamma_0
\end{equation}
Suppose $\gamma_1 < \gamma_0 < 0$, then there exists $\gamma_2$ such that $\gamma_1 < \gamma_2 < \gamma_0 < 0$. Since $\gamma_0 > \gamma_2$, for any fixed $\omega \in R$, we have

\begin{equation}
\Re \frac{g_{i_1 j_1} (j \omega) - \gamma_2 f_{i_2 j_2} (j \omega)}{f_{i_2 j_2} (j \omega)} > 0, \quad \forall (i_1 \quad j_1 \quad i_2 \quad j_2) \in I_3
\end{equation}
Namely

\begin{equation}
\frac{g_{i_1 j_1} (s) - \gamma_2 f_{i_2 j_2} (s)}{f_{i_2 j_2} (s)} \in SPR, \quad \forall (i_1 \quad j_1 \quad i_2 \quad j_2) \in I_3
\end{equation}
By Lemma 7, for any $(i_1 \quad j_1 \quad i_2 \quad j_2) \in I_3$, we have

\begin{equation}
\begin{array}{ll}
11) & \Re \frac{g_{i_1 j_1} (0) - \gamma_2 f_{i_2 j_2} (0)}{f_{i_2 j_2} (0)} > 0\\[6mm]
12) & g_{i_1 j_1} (s) - \gamma_2 f_{i_2 j_2} (s) \in H\\[6mm]
13) & f_{i_2 j_2} (s) + j \alpha [g_{i_1 j_1} (s) - \gamma_2 f_{i_2 j_2} (s)] \in H, \quad \forall \alpha \in R
\end{array}
\end{equation}
We will prove that, for any $g(s) \in K^m (s), f(s) \in K^n (s)$, we have

\begin{equation}
\frac{g(s) - \gamma_2 f(s)}{f(s)} \in SPR
\end{equation}
By Lemma 7, we only need to prove that, for any $g(s) \in K^m (s), f(s) \in K^n (s)$, we have

\begin{equation}
\begin{array}{ll}
21) & \Re \frac{g(0) - \gamma_2 f(0)}{f(0)} > 0\\[6mm]
22) & g(s) - \gamma_2 f(s) \in H\\[6mm]
23) & f(s) + j \alpha [g(s) - \gamma_2 f(s)] \in H, \quad \forall \alpha \in R
\end{array}
\end{equation}
Proof of 21): By 11) and Lemma 2, we have

\begin{equation}
f_{i_2 j_2} (0) s + g_{i_1 j_1} (0) - \gamma_2 f_{i_2 j_2} (0) \in H
\end{equation}
Namely

\begin{equation}
f_{i_2 j_2} (0) (s- \gamma_2) + g_{i_1 j_1} (0) \in H, \quad \forall (i_1 \quad j_1 \quad i_2 \quad j_2) \in I_3
\end{equation}
By Lemma 6, for any $g(s) \in K^m (s), f(s) \in K^n (s)$, we have

\begin{equation}
f(0) (s- \gamma_2) + g(0) \in H
\end{equation}
By Lemma 2, 21) is proved.

\noindent Proof of 22): By 12), we have

\begin{equation}
g_{ij} (s) - \gamma_2 f_{ij} (s) \in H, \quad \forall i,j=1,2
\end{equation}
Since $\gamma_2 < 0$, by Kharitonov Theorem\cite{ack, khar}, 22) is proved.

\noindent Proof of 23): When $\alpha = 0$, 13) becomes $f_{i_2 j_2} (s) \in H, \forall i_2, j_2 = 1, 2$. By Kharitonov Theorem\cite{ack, khar}, 23) is proved.
When $\alpha \neq 0$, 13) becomes

\begin{equation}
g_{i_1 j_1} (s) + (- \gamma_2 + \frac{1}{j \alpha}) f_{i_2 j_2} (s) \in H, \quad \forall (i_1 \quad j_1 \quad i_2 \quad j_2) \in I_3
\end{equation}
By Lemma 8, for any $g(s) \in K^m (s), f(s) \in K^n (s)$, we have

\begin{equation}
g(s) + (- \gamma_2 + \frac{1}{j \alpha}) f(s) \in H
\end{equation}
Thus, 23) is proved.
Henceforth, for any $g(s) \in K^m (s), f(s) \in K^n (s)$, we have

\begin{equation}
\frac{g(s) - \gamma_2 f(s)}{f(s)} \in SPR
\end{equation}
Namely

\begin{equation}
\Re \frac{g (j \omega)}{f (j \omega)} > \gamma_2, \quad \forall \omega \in R
\end{equation}
Hence

\begin{equation}
\inf_{\omega \in R} \Re \frac{g (j \omega)}{f (j \omega)} \geq \gamma_2
\end{equation}
Henceforth

\begin{equation}
\gamma_1 = \min \{ \inf_{\omega \in R} \Re \frac{g (j \omega)}{f (j \omega)} | g(s) \in K^m (s), f(s) \in K^n (s) \} \geq \gamma_2
\end{equation}
which contradicts $\gamma_1 < \gamma_2 < \gamma_0 < 0$. This completes the proof.

Let

\begin{equation}
\gamma_0 = \min \{ \inf_{\omega \in R} \Re \frac{g_{i_1 j_1} (j \omega)}{f_{i_2 j_2} (j \omega)} | (i_1 \quad j_1 \quad i_2 \quad j_2) \in I_3 \}
\end{equation}

\begin{equation}
\gamma_1 = \min \{ \inf_{\omega \in R} \Re \frac{g (j \omega)}{f (j \omega)} | g(s) \in K^m (s), f(s) \in K^n (s) \}
\end{equation}
Then, obviously, $\gamma_1 \leq \gamma_0$. The following result shows thay $\gamma_0$ can be regarded as an index of strict positive realness.

\vspace{0.6cm}

\noindent {\bf Theorem 4}

If $f_{ij} (s) \in H, i,j=1,2$, then

\noindent 1) \quad if $\gamma_0 < 0$, then $\gamma_1 = \gamma_0 < 0$.

\noindent 2) \quad if $\gamma_0 = 0$, then $\gamma_1 = 0$.

\noindent 3) \quad if $\gamma_0 > 0$, then $\gamma_1 \geq 0$.

Proof: 1) is a direct consequence of Theorem 3. To prove 2) and 3), suppose on the contrary that $\gamma_1 < 0$, then there exists $\gamma_2$ such that $\gamma_1 < \gamma_2 < 0$. By similar analysis as in the proof of Theorem 3, we have $\gamma_1 \geq \gamma_2$, which contradicts $\gamma_1 < \gamma_2 < 0$. This completes the proof.

\subsection{Closed Loop Systems}

Consider the open loop transfer function $\frac{g(s)}{f(s)}$. Under negative unity feedback, the closed loop transfer function is $\frac{g(s)}{f(s)+g(s)}$.

\vspace{0.6cm}

\noindent {\bf Theorem 5}

\begin{equation}
\{ \frac{g(s)}{f(s)+g(s)} | g(s) \in K^m (s), f(s) \in K^n (s) \} \subset SPR
\end{equation}
if and only if

\begin{equation}
\{ \frac{g_{i_2 j_2} (s)}{f_{i_1 j_1} (s) + g_{i_2 j_2} (s)} | (i_1 \quad j_1 \quad i_2 \quad j_2) \in I_3 \} \subset SPR
\end{equation}

Proof: Necessity: Obvious.

Sufficiency: Since $f_{ij} (s) + g_{ij} (s) \in H, \forall i,j=1,2$, by Kharitonov Theorem\cite{ack, khar}, we have

\begin{equation}
f(s) + g(s) \in H, \quad \forall g(s) \in K^m (s), f(s) \in K^n (s)
\end{equation}
By Lemma 7, we only need to prove that, for any $g(s) \in K^m (s), f(s) \in K^n (s)$, we have

\begin{equation}
\begin{array}{ll}
1) & \Re \frac{g(0)}{f(0)+g(0)} > 0\\[6mm]
2) & g(s) \in H\\[6mm]
3) & f(s) + g(s) + j \alpha g(s) \in H, \quad \forall \alpha \in R
\end{array}
\end{equation}
Proof of 1): we only need to prove

\begin{equation}
\Re \frac{f(0)+g(0)}{g(0)} > 0
\end{equation}
By Lemma 2, this is equivalent to

\begin{equation}
g(0) (s+1) + f(0) \in H
\end{equation}
By Lemma 6, this is obvious.

\noindent Proof of 2): Since $g_{ij} (s) \in H, i,j=1,2$, by Kharitonov Theorem\cite{ack, khar}, we have

\begin{equation}
g(s) \in H, \quad \forall g(s) \in K^m (s)
\end{equation}

\noindent Proof of 3): When $\alpha = 0$, 3) is obvious. When $\alpha \neq 0$, by Lemma 8, 3) is true.

\noindent This completes the proof.

From the analysis in the proof of Theorem 5, we have

\vspace{0.6cm}

\noindent {\bf Theorem 6}

For any fixed $\gamma > 0$, we have

\begin{equation}
\{ \frac{g(s)}{f(s) + \gamma g(s)} | g(s) \in K^m (s), f(s) \in K^n (s) \} \subset SPR
\end{equation}
if and only if

\begin{equation}
\{ \frac{g_{i_2 j_2} (s)}{f_{i_1 j_1} (s) + \gamma g_{i_2 j_2} (s)} | (i_1 \quad j_1 \quad i_2 \quad j_2) \in I_3 \} \subset SPR
\end{equation}

\subsection{Robust Absolute Stability}

Consider the classical Lur'e problem: the forward path is a
linear time-invariant stable component, the feedback path is a
memoryless nonlinear time-varying component (Fig. 8).

The nonlinear component $\Phi (t, \sigma)$ satisfies

\begin{equation}
\Phi (t, 0) = 0, \quad \forall t \geq 0
\end{equation}

\begin{equation}
0 \leq \sigma \Phi (t, \sigma) \leq k \sigma^2
\end{equation}
Such class of nonlinearities is said to belong to the sector $[0,
k]$, denoted as $\Phi \in sect[0, k]$. If the closed loop system
is stable for all nonlinearities in the sector $[0, k]$, then we
say the system is absolutely stable.

\vspace{0.6cm}

\noindent {\bf Lemma 9}\cite{vid}

Suppose $f(s) \in H, \Phi \in sect[0, k]$. If

\begin{equation}
\Re (\frac{1}{k} + \frac{g(j \omega)}{f(j \omega)}) > 0, \quad \forall \omega \in R
\end{equation}
Then the closed loop system is absolutely stable.

Suppose the transfer function of the forward path is an interval transfer function $\frac{g(s)}{f(s)}, g(s) \in K^m (s), f(s) \in K^n (s)$. Then we have

\vspace{0.6cm}

\noindent {\bf Theorem 7}

Suppose

\begin{equation}
\left\{ \begin{array}{ll}
        0 < k < - \frac{1}{\gamma_0} & \quad \mbox{if \, $\gamma_0 < 0$}\\[6mm]
        0 < k < + \infty             & \quad \mbox{if \, $\gamma_0 \geq 0$}
        \end{array}
\right.
\end{equation}
Then, for any $g(s) \in K^m (s), f(s) \in K^n (s), \Phi \in sect[0, k]$, the closed loop system is absolutely stable.

Proof: By Theorem 4, we have

\begin{equation}
\left\{ \begin{array}{ll}
        \gamma_1 = \gamma_0 & \quad \mbox{if \, $\gamma_0 < 0$}\\[6mm]
        \gamma_1 \geq 0     & \quad \mbox{if \, $\gamma_0 \geq 0$}
        \end{array}
\right.
\end{equation}
Hence, for any fixed $\omega \in R, g(s) \in K^m (s), f(s) \in K^n (s)$, we have

\begin{equation}
\Re (\frac{1}{k} + \frac{g(j \omega)}{f(j \omega)}) > \frac{1}{k} + \gamma_1 > 0
\end{equation}
By Lemma 9, the theorem is proved.

\section{Conclusions}

We have studied the frequency response of uncertain systems using Kharitonov stability theory on {\it first order} complex polynomial set. For an interval transfer function, we have shown that the minimal real part of the frequency response at any fixed frequency is attained at some prescribed vertex transfer functions. By further geometric and algebraic analysis, we have identified an index for strict positive realness of interval transfer functions. Some extensions and applications in positivity verification and robust absolute stability of feedback control systems have also been presented. Our results can be easily extended to the study of maximal real part, minimal (or maximal) imaginary part of the frequency response.
One salient feature of this paper is that we transform the pointwise strict positive realness problem into the stability problem of a {\it first order} complex polynomial set, thereby simplifying the analysis of the original problem. This idea is also useful in the study of maximal pointwise or bounded-bandwidth $H_{\infty}$ norm of an interval transfer function, which can be transformed into {\it Schur} stability problem of a {\it first order} complex polynomial set. For instance, ${|| \frac{g(s)}{f(s)} ||}_{\infty} < \gamma$ if and only if, for all $\omega \in R$, $f(j \omega) z + \frac{1}{\gamma} g(j \omega)$ is Schur stable. The established extreme point results will be published elsewhere.

\end{document}